%Paper: dg-ga/9503002
%From: Dan Burghelea <burghele@math.ohio-state.edu>
%Date: Wed, 8 Mar 1995 14:26:24 -0500

\input amstex
\documentstyle{amsppt}
\magnification=1200
\topmatter
\title Mayer-Vietoris Formula for Determinants of Elliptic Operators of
Laplace-Beltrami Type.
(after Burghelea, Friedlander and Kappeler)
\endtitle
\rightheadtext{Mayer-Vietoris Formula for Determinants of Elliptic Operators}
\author Yoonweon Lee
\endauthor
\date March 5, 1995
\enddate
\abstract
The purpose of this note is to provide a short cut presentation of a
Mayer-Vietoris formula due to Burghelea-Friedlander-Kappeler for
the regularized determinant in the case of elliptic operators of
Laplace Beltrami type in the form typically needed in applications to torsion.

\endabstract
\endtopmatter
\def\dd{\vskip 0.5 cm}
\def\ddd{\vskip 0.2 cm}
\def\dddd{\vskip 0 cm}
\document

{\bf 1. Statement of Mayer-Vietoris Formula for Determinants}
\ddd
Let $(M,g)$ be a closed oriented Riemannian manifold of dimension $d$ and
$\Gamma$ be an oriented submanifold of codimension 1. We denote by $\nu$
the unit normal vector field along $\Gamma$. Let $M_{\Gamma}$
be the compact manifold with boundary $\Gamma^{+}\sqcup\Gamma^{-}$ obtained
by cutting $M$ along $\Gamma$, where $\Gamma^{+}$ and $\Gamma^{-}$ are copies
of $\Gamma$
and denote by $p:M_{\Gamma}\rightarrow M$ the identification map.
 The vector field $\nu$ has the lift on $M_{\Gamma}$ which we denote by $\nu$
again.
Denote by $\Gamma^{+}$
the component of the boundary where the lift of $\nu$ points outward.
Given a smooth vector bundle $E\rightarrow M,$ denote by $E_{\Gamma}$
the pull back of $E\rightarrow M$ to $M_{\Gamma}$ by $p$.
Let $A:C^{\infty}(E)\rightarrow C^{\infty}(E)$ be an elliptic,
essentially self-adjoint,
positive definite differential operator of Laplace-Beltrami type,
where we say that $A$ is of Laplace-Beltrami type if $A$
is an operator of order 2
whose principal symbol is $\sigma_{L}(x,\xi)=\parallel\xi\parallel^{2}Id_x$,
$Id_x\in End_x(E_x,E_x)$.
We denote by $A_{\Gamma}:C^{\infty}(E_{\Gamma})\rightarrow
C^{\infty}(E_{\Gamma})$ the extension
of $A$ to smooth sections of $E_{\Gamma}$.

Consider Dirichlet and Neumann boundary conditions $B,C$ on
$\Gamma^{+}\sqcup\Gamma^{-}$ defined as follows:
$$B:C^{\infty}(E_{\Gamma})\rightarrow C^{\infty}(E_{\Gamma}
\mid_{\Gamma^{+}\sqcup\Gamma^{-}}), B(f)=f\mid_{\Gamma^{+}\sqcup\Gamma^{-}}$$
$$C:C^{\infty}(E_{\Gamma})\rightarrow C^{\infty}(E_{\Gamma}
\mid_{\Gamma^{+}\sqcup\Gamma^{-}}), C(f)=\nu(f)\mid_{\Gamma^{+}\sqcup
\Gamma^{-}}.$$

Consider $\tilde {A_{\Gamma,B}}=(A_{\Gamma},B):C^{\infty}(E_{\Gamma})
\rightarrow
C^{\infty}(E_{\Gamma})\oplus
C^{\infty}(E_{\Gamma}\mid_{\Gamma^{+}\sqcup\Gamma^{-}})$.
{}From the properties of $A$ it follows that
 $\tilde{A_{\Gamma,B}}$ is invertible. Therefore we can define the
corresponding Poisson operator
$P_{B}$ as the restriction of $\tilde
{A_{\Gamma,B}}^{-1}$ to $0\oplus
C^{\infty}(E_{\Gamma}\mid_{\Gamma^{+}\sqcup\Gamma^{-}})$.
Denote by $A_{B}$ the
restriction of $\tilde {A_{\Gamma,B}}$ on $\{u\in C^{\infty}(E_{\Gamma})\mid
B(u)=0\}$. Then $A_{B}$ is also essentially self-adjoint
and positive definite (cf Lemma 3.1). This ( using  standard
analytic continuation technique due to Seeley (cf.[Se])) allows us to define
$$log Det(A)= -\frac{d}{ds}\mid_{s=0}tr\frac{1}{2\pi i}
\int_{\gamma}\lambda^{-s}(\lambda-A)^{-1}d\lambda$$
$$logDet(A_{\Gamma},B)= -\frac{d}{ds}\mid_{s=0}tr\frac{1}{2\pi i}\int_{\gamma}
\lambda^{-s}(\lambda-A_{B})^{-1}d\lambda,$$
where $\gamma$ is a path around the negative real axis,
$$\{\rho e^{i\pi}\mid \infty >\rho\geq \epsilon\}\cup
\{\epsilon e^{i\theta}\mid \pi\geq \theta\geq -\pi\}\cup \{\rho e^{-i\pi}\mid
\epsilon\leq\rho <\infty\}$$ with $\epsilon >0$ chosen sufficiently small to
ensure that $\Gamma$ does not separate the spectrum.

We define the Dirichlet to Neumann operator, associated to
A,B and C,
$$R:C^{\infty}(E\mid _{\Gamma})
\rightarrow C^{\infty}(E\mid_{\Gamma})$$ by the composition of the following
maps
$$C^{\infty}(E\mid_{\Gamma})@>\bigtriangleup_{ia}>>
C^{\infty}(E\mid_{\Gamma^{+}})\oplus C^{\infty}(E\mid_{\Gamma^{-}})@>P_{B}>>
 C^{\infty}(E_{\Gamma})@>C >> C^{\infty}
(E\mid_{\Gamma^{+}})\oplus C^{\infty}(E\mid_{\Gamma^{-}})$$
$$@>\bigtriangleup_{if} >> C^{\infty}(E\mid_{\Gamma}),$$
where $\bigtriangleup_{ia}(f)=(f,f)$ is the diagonal inclusion and
$\bigtriangleup_{if}(f,g)= f - g$ is the difference operator.
Then $R$ is an essentially self-adjoint, positive definite,
elliptic operator of order 1
(cf Lemma 3.5).
\proclaim{Theorem 1.1 (Mayer-Vietoris Type Formula for Determinants [BFK])}

Let $(M,g)$ be a closed oriented Riemannian manifold of dimension $d$ and
$A$ be an elliptic, essentially self-adjoint, positive definite differential
operator
of Laplace-Beltrami type acting on smooth sections of a vector bundle
$E\rightarrow M$.
Then $A_{B}$ and $R$ are
essentially self-adjoint, positive definite elliptic operators and
$$Det(A)=cDet(A_{\Gamma},B)Det(R),$$
where $c$ is a local quantity which can be computed in terms of the symbols
of $A,B$ and $C$ along $\Gamma$.
\endproclaim

Remark: The above result can be extended to manifolds with boundary. E.g.
consider an oriented, compact, smooth manifold $M$ whose boundary
$\partial M$ is a disjoint union of two components $\partial_{+} M$
and $\partial_{-} M$ with $\Gamma \cap \partial M = \emptyset,$
an operator $A$ of Laplace-Beltrami type and
differential elliptic boundary conditions $B_{+}$ respectively
$B_{-}$ for $A$ on $\partial_{+} M$
respectivley $\partial_{-} M.$ Denote by $A^{(0)}$ the operator $A$ with
domain $\{ u \in C^{\infty}(E) \mid B_{+}u=0, B_{-}u=0 \}.$ Then Theorem 1.1
remains true with $A$ replaced by $A^{(0)}.$

\dd
{\bf 2. The Asymptotics of Determinants of Elliptic Pseudodifferential
Operators with Parameter}
\ddd
 Let V be an open angle in the complex $\lambda$-plane and P($\lambda$),
$\lambda\in $V, a family of $\Psi$DO's of order m, m a positive integer,
acting on smooth sections of
 a vector bundle $E\rightarrow M$ of rank $\nu$, where $M$ denotes a closed
smooth Riemannian manifold of dimension $d$.

\proclaim{Definition 2.1}(cf.[Sh]) The family
 P($\lambda$), $\lambda \in V,$ is said to be a $\Psi$DO with parameter of
weight $\chi >0$ if in any coordinate neighborhood U of $M$, not necessarily
connected, and for an arbitarily fixed $\lambda \in V$, the complete symbol
$p(\lambda;x,\xi)$ of P is in
 $C^{\infty}(U\times \Bbb R^{d},End(\Bbb C^{\nu}))$
and, moreover, for any  multiindices
$\alpha$ and $\beta$, there exists a
constant $C_{\alpha,\beta}$ such that
$\mid \partial^{\alpha}_{\xi}\partial^{\beta}_{x}p(\lambda;x,\xi)\mid\leq C
_{\alpha\beta}(1+\mid\xi\mid+\mid\lambda
\mid^{\frac{1}{\chi}})^{m-\mid\alpha\mid}$.
\endproclaim
\proclaim{Definition 2.2}
 P($\lambda$) is called classical if in any chart the complete symbol \newline
$p(\lambda;x,\xi)$ admits an expansion of the form
$$ p(\lambda;x,\xi)\sim p_{m}(\lambda;x,\xi)+p_{m-1}(\lambda;x,\xi)+\cdots,$$
where $p_{j}(\tau^{\chi}\lambda;x,\tau\xi)=\tau^{j}p_{j}(\lambda;x,\xi)
 (\tau>0,j\leq m)$.
The family P($\lambda$) is said to be elliptic with parameter if
$p_{m}(\lambda;x,\xi)$  is invertible for all $x\in M$, $\xi\in
T^{\ast}_{x}(M)$ and $\lambda\in$V  satisfying
$\mid\xi\mid + \mid\lambda\mid^{\frac{1}{\chi}}\neq 0$.
\endproclaim

\proclaim{Definition 2.3}
Let $Q$ be an elliptic $\Psi$DO. The angle $\pi$ is called an Agmon angle
for $Q$ if for some $\epsilon>0$, $spec(Q)\cap\Lambda_{\epsilon}=\emptyset$,
where $spec(Q)$ denotes the spectrum of $Q$ and $\Lambda_{\epsilon}
=\{z\in\Bbb C\mid \pi-\epsilon<arg(z)<\pi+\epsilon \text{ or }
\mid z\mid<\epsilon\}$.
\endproclaim

\proclaim{Theorem 2.4}
 Let $P(\lambda)$ be an essentially self-adjoint, positive definite,
 classical $\Psi$DO of
order $m\in\Bbb N$ with parameter
$\lambda\in V$ of weight $\chi >0$ such that \newline
{\rm{(i)}} $P(\lambda)$ is elliptic with parameter and \newline
{\rm{(ii)}} for each $\lambda\in V$, $P(\lambda)$ has $\pi$ as an Agmon angle.
\newline
 Then $logDetP(\lambda)$ admits an asymptotic expansion for
$\lambda\in V$, $\mid\lambda\mid\rightarrow\infty$, of the form
$$logDetP(\lambda)\sim \sum^{\infty}_{j=-d}\pi_{j}
\mid\lambda\mid^{-\frac{j}{\chi}} +
\sum^{d}_{j=0}q_{j}\mid\lambda\mid^{\frac{j}{\chi}}log\mid\lambda\mid.$$
 The coefficients $\pi_{j}$ and $q_{j}$ can be evaluated in terms of the symbol
of $P$ and $\frac{\lambda}{\mid\lambda\mid}$. In particular, $\pi_{0}$
is independent of perturbations by lower order operators, whose orders differ
at least by $d+1$ from the order of $P(\lambda)$.
\endproclaim

For the convenience of the reader we include the proof of this theorem
which can be found in the appendix of [BFK].
\ddd
{\bf Proof of Theorem 2.4}
 We divide the proof into several steps.

{\bf Step 1}
 By a standard procedure we construct a parametrix for \newline
$R(\mu,\lambda)=(\mu -P(\lambda))^{-1} (\mu\leq 0)$.

{\bf Step 2}
 Define $R_{N}(\mu,\lambda)$ to be a conveniently chosen approximation of
$R(\mu,\lambda)$ and write $P(\lambda)^{-s}=
P_{N}(\lambda;s)+\tilde{P_{N}}(\lambda;s)$, where $P_{N}(\lambda;s)=
\frac{1}{2\pi i}\int_{\gamma}\mu^{-s}R_{N}(\mu,\lambda)d\mu$ and where
$\gamma$ denotes a contour around the negative axis, enclosing the origin in
clockwise orientation. Then for $s$ in $\Bbb C$ with $Re s$ sufficiently
large, $\zeta(s)=\zeta_{N}(s)+\tilde{\zeta_{N}}(s)$,
where $\zeta(s)=trP(\lambda)^{-s}, \zeta_{N}(s)=trP_{N}(\lambda,s)$, and
$\tilde{\zeta_{N}}(s)=tr\tilde{P_{N}}(\lambda,s)$.

{\bf Step 3}
 Describe an asymptotic expansion of
$\frac{\partial}{\partial s}\mid _{s=0}\zeta_{N}
(\lambda,s)$ as $\lambda\rightarrow \infty$.

{\bf Step 4}
 Provide an estimate for the remainder term
$\frac{\partial}{\partial s}\mid _{s=0}
\tilde{\zeta _{N}}(\lambda,s)$ as $\lambda\rightarrow \infty$.

{\bf Step 5}
Provide a formula for $\pi _{0}.$

\ddd
{\bf Step 1}
 We want to construct a parametrix for $R(\mu,\lambda)=(\mu-P(\lambda))^{-1}$
\newline
$(\mu\leq 0)$. Consider the equation
$(\mu -p(\lambda;x,\xi))\circ r(\mu,\lambda;x,\xi)=Id$, where
$\circ$ denotes mulplication in the algebra of symbols.

 Introduce, for $\alpha =(\alpha_{1},\cdots,\alpha_{d})$, the standard notation
$\alpha !=\alpha_{1}!\cdots\alpha_{d}!$, $\partial^{\alpha}_{\xi}=
(\frac{\partial}{\partial\xi_{1}})^{\alpha_{1}}\cdots
(\frac{\partial}{\partial\xi_{d}})^{\alpha_{d}}$ and $D^{\alpha}_{\xi}=
(\frac{1}{i})^{\alpha}\partial^{\alpha}{\xi}$.

 Write $r(\mu,\lambda;x,\xi)\sim r_{-m}(\mu,\lambda;x,\xi) +
r_{-m-1}(\mu,\lambda;x,\xi) +\cdots$, where $r_{-j}(\mu,\lambda;x,\xi)$ is
positive homogeneous of degree
$-j$ in $(\xi,\mu ^{\frac{1}{m}},\lambda ^{\frac{1}{\chi}})$. Then we obtain
the following formula:
$$r_{-m}(\mu,\lambda;x,\xi)= (\mu - p_{m}(\lambda;x,\xi))^{-1}$$
and for $j\geq 1$,
$$r_{-m-j}(\mu,\lambda;x,\xi)$$
$$= -(\mu - p_{m}(\lambda;x,\xi))^{-1}\sum^{j-1}_{k=0}
\sum_{\mid\alpha\mid +l+k=j}\frac{1}{\alpha !}\partial ^{\alpha}_{\xi}
p_{m-l}(\lambda;x,\xi)D^{\alpha}_{x}r_{-m-k}(\mu,\lambda;x,\xi).$$
 The functions $r_{j}(\mu,\lambda;x,\xi)$ satisfy the following homogeneity
condition
$$r_{j}(\tau ^{m}\mu,\tau ^{\chi}\lambda;x,\tau\xi)=
\tau ^{j}r_{j}(\mu,\lambda;x,\xi) (\tau >0).$$

 By a standard procedure, $\sum_{j\geq 0}r_{-m-j}(\mu,\lambda;x,\xi)$ gives
rise to a $\Psi$DO with parameter, called a parametrix for $R(\mu,\lambda)$.

{\bf Step 2}
 Introduce a finite cover ($U_{j}$) of $M$ by open charts and take a
partition of unity $\varphi_{j}$, subordinate to $U_{j}$.
Choose $\psi_{j}\in C^{\infty}_{0}(U_{j})$ such that $\psi_{j}\equiv 1$
in some neighborhood of supp$\varphi_{j}$. Let us fix local coordinates in
every $U_{j}$ and define the operators
$$(R^{(j)}_{N,\mu,\lambda}f)(x)=
\psi_{j}(x)\cdot \int_{\Bbb R^{d}}d\xi \frac{1}{(2\pi)^{d}}\int_{\Bbb R^{d}}dy
\left(r^{(N)}(\mu,\lambda;x,\xi)e^{i(x-y)\cdot\xi}\varphi_{j}(y)f(y)\right),$$
where $r^{(N)}(\mu,\lambda;x,\xi)= \sum^{N-1}_{j=0}r_{-m-j}(\mu,\lambda;x,\xi)$
in the local coordinates of $U_{j}$. The approximation $R_{N}(\mu,\lambda)$ of
the resolvent $R(\mu,\lambda)$ is defined by $R_{N}(\mu,\lambda)=$ \newline
$\sum_{j}R^{(j)}_{N}(\mu,\lambda)$.
We need an estimate of $R(\mu,\lambda)-R_{N}(\mu,\lambda)$ in trace norm.
The latter is denoted by $\mid \parallel\cdot\parallel \mid$.

\proclaim{Lemma 2.5}
 Choose $N > \frac{3d}{2}+ m$. Then for $\lambda\in V_{1}$ and $\mu\in
\Bbb R^{-}$
with $\mid \mu\mid$ sufficiently large
$$\mid \parallel R(\mu,\lambda)- R_{N}(\mu,\lambda)\parallel \mid <
C_{N}(1+\mid\lambda\mid)^{-\frac{(N-\frac{3d}{2}-m)}{\chi}}
(1+\mid\lambda\mid)^{-2},$$
where $V_{1}$ is an angle whose closure is contained in $V$, $V_{1}\ll V.$
\endproclaim
\demo{ Proof}
 Define $T_{N}(\mu,\lambda)$ by $(\mu- P(\lambda))R_{N}(\mu,\lambda)=
Id- T_{N}(\mu,\lambda)$. From $(\mu- P(\lambda))\cdot$ \newline
$R(\mu,\lambda)=Id$, we then conclude that
$R(\mu,\lambda)- R_{N}(\mu,\lambda)= R(\mu,\lambda)T_{N}(\mu,\lambda)$.

 The claimed estimate of the lemma follows, once we have proved that for some
$\tau > d$
$$\parallel R(\mu,\lambda)\parallel_{L^{2}\rightarrow
L^{2}}\leq C(1+\mid\mu\mid)^{-1} \eqno (2.1)$$
$(\lambda\in V_{1}\ll V,\mu\in \Bbb R^{-},\mid\mu\mid$ sufficiently large)
and  \newline
$$\parallel T_{N}(\mu,\lambda)\parallel _{L^{2}\rightarrow
H^{\tau}}\leq C_{\tau}(1+\mid\lambda\mid^{\frac{1}{\chi}}+
\mid\mu\mid^{\frac{1}{m}})^{-N+\tau} \eqno (2.2)$$
 because, from (2.2) we can conclude that $T_{N}(\mu,\lambda)$ is a
$\Psi$DO of order $-\tau < -d$ and hence of trace class, when considered as
an operator on $L^{2}$-sections of $E\rightarrow M$.
The estimate (2.1) is standard (cf.e.g.[Sh]) and (2.2) follows from the fact
that the symbol $t_{N}(\mu,\lambda;x,\xi)$ of $T_{N}(\mu,\lambda)$ satisfies
$$\mid D^{\alpha}_{x}D^{\beta}_{\xi}t_{N}(\mu,\lambda;x,\xi)\mid
\leq C_{\alpha\beta}(1+\mid\xi\mid + \mid\lambda\mid^{\frac{1}{\chi}}+
\mid\mu\mid^{\frac{1}{m}})^{-N-\mid\beta\mid}.$$
Thus the norm of $T_{N}(\mu,\lambda)$ as an operator from $H^{s}$ to
$H^{s+\tau}$ is $$O(1+\mid\lambda\mid^{\frac{1}{\chi}}+
\mid\mu\mid^{\frac{1}{m}})^{-N+\tau}.$$ Therefore
$$\mid\parallel T_{N}\parallel\mid = O(1+\mid\lambda\mid^{\frac{1}{\chi}}+
\mid\mu\mid^{\frac{1}{m}})^{-N+\frac{3d}{2}}.$$
\qed
\enddemo

{\bf Step 3}
 Next we study the asymptotic expansion of
$\frac{\partial}{\partial s}\mid _{s=0}\zeta_{N}(\lambda,s)$ as
$\lambda\rightarrow +\infty$.
Recall that $P_{N}(\lambda,s)=
\frac{1}{2\pi i}\int_{\gamma}d\mu \mu^{-s}R_{N}(\mu,\lambda)$.
Its Schwarz kernel is given by
$$ P_{N}(\lambda,s;x,y)=\sum_{j}\psi_{j}(x)\varphi_{j}(y)
\frac{1}{(2\pi)^{d}}\int_{\Bbb R^{d}}d\xi e^{i(x-y)\cdot\xi}
\frac{1}{2\pi i}\int_{\gamma}d\mu \mu^{-s}\sum^{N-1}_{k=0}r_{-m-k}.$$
As a consequence
$$P_{N}(\lambda,s;x,x)=
\sum_{j}\varphi_{j}(x)\frac{1}{(2\pi)^{d}}\sum^{N-1}_{k=0}
I_{k}(s,\lambda,x)=\frac{1}{(2\pi)^{d}}\sum^{N-1}_{k=0}I_{k}(s,\lambda,x),$$
where $I_{j}(s,\lambda,x)=\frac{1}{2\pi i}\int_{\Bbb R^{d}}d\xi
\int_{\gamma}d\mu \mu^{-s}r_{-m-j}(\mu,\lambda;x,\xi)$.
By the change of variables $\xi =\mid\lambda\mid^{\frac{1}{\chi}}\xi',
\mu=\mid\lambda\mid^{\frac{m}{\chi}}\mu'$ and by using the homogeneity of
$r_{-m-j}$, we obtain
$$I_{j}(s,\lambda,x)=\frac{1}{2\pi i}\mid\lambda\mid^{\frac{d-ms-j}{\chi}}
\int_{\Bbb R^{d}}d\xi\int_{\gamma}d\mu \mu^{-s}r_{-m-j}
(\mu,\frac{\lambda}{\mid\lambda\mid};x,\xi).$$
 We need to investigate $J_{k}(s,\lambda,x):=\frac{1}{2\pi i}
\int_{\Bbb R^{d}}d\xi\int_{\gamma}d\mu \mu^{-s}r_{-m-k}
(\mu,\frac{\lambda}{\mid\lambda\mid};x,\xi)$.

\proclaim{Lemma 2.6}
 Let $\omega\in V$ with $\mid\omega\mid =1$. Then $J_{k}(s,\omega;x)$ is
holomorphic in $s$ in the half plane $Re s >\frac{d-k}{m}$ and it admits a
meromorphic continuation in the complex s-plane. The point s=0 is always
regular and $J_{k}(0,\omega;x)=0$ if $k > d$.
\endproclaim
\demo{ Proof}
(i) Integrating by parts with respect to $\mu$, one obtains
$$J_{k}(s,\omega;x)=\int_{\Bbb R^{d}}d\xi\frac{1}{(1-s)\cdots(l-s)}
\frac{1}{2\pi i}\int_{\gamma}(-1)^{l}d\mu \mu^{-s+l}
\frac{\partial^{l}}{\partial\mu^{l}}r_{-m-k}(\mu,\omega;x,\xi).$$
If $l> Res-1$, the contour integral reduces to
$$-\frac{sin\pi s}{\pi}
\int^{\infty}_{0}d\mu \mu^{-s+l}(\frac{\partial^{l}}{\partial\mu^{l}}
r_{-m-k})(-\mu,\omega;x,\xi).$$
Further, the matrix $\frac{\partial^{l}}{\partial\mu^{l}}r_{-m-k}$ can be
estimated
$$\mid\frac{\partial^{l}}{\partial\mu^{l}}r_{-m-k}\mid\leq
C(1+\mid\mu\mid^{\frac{1}{m}}+\mid\xi\mid)^{-m-k-ml}.$$
Thus, for $Res > \frac{d-k}{m}$, the integral
$$\int_{\Bbb R^{d}}d\xi\int^{\infty}_{0}d\mu \mu^{-s+l}
((\frac{\partial^{l}}{\partial\mu^{l}})r_{-m-k})(-\mu,\omega;x,\xi)$$
converges absolutely and therefore is a holomorphic function in s. Moreover,
$-\frac{sin\pi s}{\pi (1-s)\cdots (l-s)}$  is entire. In all, we have proved
that $J_{k}(s,\lambda;x)$ is holomorphic in $Res > \frac{d-k}{m}$.

(ii) Next let us prove that $J_{k}(s,\lambda;x)$ can be meromorphically
continued to the entire complex $s$-plane. To keep the exposition simple let
us assume that $P(\lambda)$ is a scalar $\Psi$DO. The expressions
$r_{-m-k}(\mu,\omega;x,\xi)$ have been defined in a recursive fashion
and are sums of terms of the form
$(\mu-p_{m}(\omega;x,\xi))^{-l}q_{l,k}(\omega;x,\xi)$ with $l\geq 1$, where
ord$(q_{l,k})=-m-k+ml$ and $q_{l,k}$ is an expression, independent of
$\mu$, involving only the symbols $p_{m-j}(\omega;x,\xi)$ and their
derivatives with $0\leq j\leq k$.

 It follows from the recursive definition of the $r_{-m-k}$ that $l$ has to
satisfy $l\geq k+1$ and thus, in the case $k\geq 1$, $J_{k}$ consists of a
sum of terms of the form
$$\int_{\Bbb R^{d}}d\xi q_{l,k}(\omega;x,\xi)\frac{1}{2\pi i}
\int_{\gamma}d\mu \mu^{-s}(\mu-p_{m}(\omega;x,\xi))^{-l}$$
$$=\left(\int_{\Bbb R^{d}}d\xi q_{l,k}(\omega;x,\xi)(p_{m}
(\omega;x,\xi))^{-s-l+1}\right)\left(\frac{(-1)^{l-1}}{(l-1)!}s(s+1)\cdots
(s+l-2)\right),$$
where after integration by parts, we used Cauchy's formula.
 As $\mid\omega\mid=1$, it follows from Definition 2.1 that the integrand
$q_{l,k}(\omega;x,\xi)p_{m}(\omega;x,\xi)^{-s-l+1}$ is absolutely integrable
in $\mid\xi\mid\leq 1$. Thus one only needs to consider the integral over
$\mid\xi\mid >1$. For $\omega$ fixed, the symbols $q_{l,k}$ and $p_{m}$ are
classical and admit an asymptotic expansion in $\xi$-homogeneous functions.
Consider two cases:

Case 1: $k=0$.

 Using that $r_{-m}(\mu,\omega;x,\xi)=(\mu-p_{m}(\omega;x,\xi))^{-1}$ we
conclude that
$$J_{0}(s,\omega;x)=\int_{\Bbb R^{d}}d\xi\frac{1}{2\pi i}
\int_{\gamma}d\mu \mu^{-s}(\mu-p_{m}(\omega;x,\xi))^{-1}=
\int_{\Bbb R^{d}}d\xi(p_{m}(\omega;x,\xi))^{-s}.$$
 Recall that $\omega$ with $\mid\omega\mid=1$ is fixed and thus
$p_{m}(\omega;x,\xi)$
defines an elliptic $\Psi$DO $P_{m}(\omega;x,D)$ and we can apply the
standard theory of complex powers of elliptic operators (cf.e.g.[Se]) to
conclude that $\int_{\Bbb R^{d}}d\xi p_{m}(\omega;x,\xi)^{-s}$ has a
meromorphic continuation in the whole complex $s$-plane, with at most simple
poles and that $s=0$ is a regular point. The poles are located at
$s_{j}=\frac{d-j}{m}$ with $j\in\{0,1,2,\cdots\} \setminus
\{d\}$.

Case 2: $k\geq 1$.

 As it was observed by Guillemin [Gu] and Wodzicki [Wo] in the context of
non-commutative residues, $\int_{\Bbb R^{d}}q_{l,k}(\omega;x,\xi)(p_{m}
(\omega;x,\xi))^{-s-l+1}d\xi$ admits a meromorphic continuation to the whole
complex $s$-plane with
at most simple poles. Thus \newline
$s\cdot\int_{\Bbb R^{d}}d\xi q_{l,k}
(\omega;x,\xi)(p_{m}(\omega;x,\xi))^{-s-l+1}$ must be regular at $s=0$.
This shows that $J_{k}(s,\omega;x)(k\geq 1)$ is meromorphic and that $s=0$ is
a regular point.

(iii) Let $k > d+1-m$. Observe that
$$\mid r_{-m-k}(\mu,\omega;x,\xi)\mid\leq C_{k}(1+\mid\mu\mid^{\frac{1}{m}}+
\mid\xi\mid)^{-m-k}.$$
As $m\geq 1$, the integral
$J_{k}(s,\omega;x)=
\frac{1}{2\pi i}\int_{\Bbb R^{d}}d\xi\int_{\gamma}d\mu \mu^{-s}
r_{-m-k}(\mu,\omega;x,\xi)$ converges absolutely at $s=0$. Evaluating at
$s=0$, one obtains $\int_{\gamma}d\mu r_{-m-k}(\mu,\omega;x,\xi)=0$ \newline
and thus $J_{k}(0,\omega;x)=0$ for $k>d$.
\qed
\enddemo

By the above lemma, we see that
$$P_{N}(\lambda,s;x,x)=\frac{1}{(2\pi)^{d}}\sum^{N-1}_{k=0}\mid\lambda\mid
^{\frac{(d-ms-k)}{\chi}}J_{k}(s,\frac{\lambda}{\mid\lambda\mid};x).$$
Hence, with $N^{\ast}=min(N-1,d)$,
$$\frac{\partial}{\partial s}\mid _{s=0}P_{N}(\lambda,s;x,x)=
\frac{1}{(2\pi)^{d}}\sum^{N-1}_{k=0}\mid\lambda\mid^{\frac{d-k}{\chi}}
\frac{\partial}{\partial s}
J_{k}(s,\frac{\lambda}{\mid\lambda\mid};x)\mid_{s=0}-$$
$$\frac{m}{\chi}\frac{1}{(2\pi)^{d}}
\sum^{N^{\ast}}_{k=0}\mid\lambda\mid^{\frac{d-k}{\chi}}
log\mid\lambda\mid\cdot J_{k}(0,\frac{\lambda}{\mid\lambda\mid};x).$$

{\bf Step 4}
 We have to estimate $tr\tilde{P_{N}}(\lambda,s)$,
where
$$\tilde{P_{N}}(\lambda,s)= P(\lambda)^{-s}-P_{N}(\lambda,s)=
\frac{1}{2\pi i}\int_{\gamma}d\mu \mu^{-s}(R(\mu,\lambda)-R_{N}(\mu,\lambda)).
$$

The estimate of Lemma 2.5 implies that
$$\mid tr\frac{\partial}{\partial s}\tilde{P_{N}}(\lambda,s)\mid_{s=0}
\leq C(1+\mid\lambda\mid)^{-\frac{(N-\frac{3}{2}d-m)}{\chi}}
\int^{\infty}_{0}d\mu\cdot \frac{\mid log\mu\mid}{1+\mid\mu\mid^{2}}$$
$$\leq C(1+\mid\lambda\mid)^{-\frac{(N-\frac{3}{2}d-m)}{\chi}}$$
and thus the asymptotic expansion given in Theorem 2.4 is proved.

{\bf Step 5}
 In the notation introduced above, we obtain the following formula for
$\pi_{0}$:
$$\pi_{0}=\sum_{j}\frac{\partial}{\partial s}\frac{1}{(2\pi)^{d}}
\int_{\Bbb R^{d}}dvol(x)J_{d}(s,\lambda;x)\varphi_{j}(x)\mid _{s=0},$$
where $J_{d}(s,\lambda;x)=\frac{1}{2\pi i}\int_{\Bbb R^{d}}d\xi
\int_{\gamma}d\mu \mu^{-s}r_{-m-d}
(\mu,\frac{\lambda}{\mid\lambda\mid};x,\xi)$.

As $r_{-m-d}(\mu,\frac{\lambda}{\mid\lambda\mid};x,\xi)$ is defined recursively
by $(j\geq 1)$
$$r_{-m-j}(\mu,\lambda;x,\xi)=$$
$$-(\mu-p_{m}(\lambda;x,\xi))^{-1}\sum^{j-1}_{k=0}
\sum_{\mid\alpha\mid +l+k=j}\frac{1}{\alpha !}\partial^{\alpha}_{\xi}p_{m-l}
(\lambda;x,\xi)D^{\alpha}_{x}r_{-m-k}(\mu,\lambda;x,\xi),$$
we conclude that $\pi_{0}$ only depends on
$p_{m-j}(\frac{\lambda}{\mid\lambda\mid};x,\xi)$ for
$0\leq j\leq d$ and its derivatives up to order $d$.
\qed

The following result is due to Voros [Vo] and Friedlander [Fr].
For the convenience of the reader we include Voros' proof.
\proclaim{Proposition 2.7}
Let $\{\lambda_{k}\}_{k\geq 1}$ be a sequence in $V_{0,\frac{\pi}{2}}=
\{z\in\Bbb C\mid -\frac{\pi}{2}+\epsilon<arg(z)<\frac{\pi}{2}-\epsilon\},$
 possibly with multiplicities,
arranged in such a way that $0<Re \lambda_{1}\leq Re \lambda_{2}\leq \cdots.$
Assume that the heat trace, $\theta(t):=\sum^{\infty}_{k=0}e^{-t\lambda_{k}}
(t>0)$ admits an asymptotic expansion for $t\rightarrow 0$ of the form
$$\theta(t)\sim\sum_{n\geq 0}c_{i_{n}}t^{i_{n}},$$
where $i_{0}<0$ and $i_{0}<i_{1}<i_{2}<\cdots\rightarrow +\infty$.
For $\lambda$ with $Re \lambda>0$, let $\zeta(s,\lambda)=\sum^{\infty}_{k=0}
(\lambda_{k}+\lambda)^{-s}$.
Then $\pi_{0}=0$, where $\pi_{0}$ is the constant term in the asymptotic
expansion of $-\frac{d}{ds}\zeta(s,\lambda)\mid_{s=0}$ for $\mid\lambda\mid
\rightarrow +\infty$.
\endproclaim
\demo{Proof}
It is well known that for $Re s>-i_{0}$, $\zeta(s,\lambda)$ is well defined and
$\zeta(s,0)=
\frac{1}{\Gamma(s)}\int^{\infty}_{0}\theta(t)\cdot t^{s-1}dt.$ Let
$$\eta(s,\lambda)=\int^{\infty}_{0}\sum_{k=0}^{\infty}e^{-\lambda_{k}t}
e^{-\lambda t}t^{s-1}dt.$$
Then
$\zeta(s,\lambda)=\frac{1}{\Gamma(s)}\eta(s,\lambda)$.
For $Re s>-i_{0}$, $\eta(s,\lambda)$ can be expanded in $\lambda$ for
$\mid\lambda\mid\rightarrow\infty$
$$\eta(s,\lambda)\sim\sum_{n\geq 0}c_{i_{n}}\lambda^{-s-i_{n}}\int_{0}^{\infty}
t^{i_{n}+s-1}e^{-t}dt=\sum_{n\geq 0}c_{i_{n}}\Gamma(s+i_{n})
\lambda^{-s-i_{n}}.$$
Thus
$$\zeta(s,\lambda)=\frac{1}{\Gamma(s)}\eta(s,\lambda)
\sim\frac{1}{\Gamma(s)}\lambda^{-s}\sum_{n\geq 0}c_{i_{n}}\Gamma(s+i_{n})
\lambda^{-i_{n}}.$$
All functions involved are meromorphic functions of $s$. Moreover $s=0$ is a
regular point of $\zeta(s,\lambda)$ and thus
$\frac{d}{ds}\zeta(s,\lambda)\mid_{s=0}$ admits an asymptotic expansion in
$\lambda$ of the form
$$\frac{d}{ds}\zeta(s,\lambda)\mid_{s=0}
\sim \sum_{i_{n}\notin \Bbb Z^{-}\cup\{0\}}c_{i_{n}}
\Gamma(i_{n})\lambda^{-i_{n}}+\sum_{i_{n}\in\Bbb Z^{-}}c_{i_{n}}\frac{d}{ds}
\frac{1}{(s+i_{n})\cdots(s-1)}\mid_{s=0}\lambda^{-i_{n}}-$$
$$\sum_{i_{n}\in\Bbb Z^{-}\cup\{0\}}\frac{c_{i_{n}}}{(i_{n})\cdots(-1)}
\lambda^{-i_{n}}log\lambda,$$
where $\Bbb Z^{-}$ is the set of all negative integers.
This expansion shows that $\pi_{0}=0$.
\qed
\enddemo

\dd
{\bf 3. Auxiliary Results for the Proof of Theorem 1.1}
\ddd
We begin by collecting a number of results about operators related to
$A$ and $\Gamma$. Denote by $H^s(E_{\Gamma})$ the Sobolev spaces of
$E_{\Gamma}-$valued sections. Throughout section 3 and section 4
we assume that $A$ satisfies the hypothesis of Theorem 1.1 and
fix $\epsilon > 0,$ so that the spectrum of $A$ is bounded from
below by $\epsilon.$
\proclaim{Lemma 3.1}
{\rm{(i)}} The operator
$A_{B}:\{u\in C^{\infty}(E_{\Gamma})\mid B(u)=0\}\rightarrow C^{\infty}
(E_{\Gamma})$
has a  self-adjoint extension $\bar{A}_{B}$ with domain
 $D(\bar{A}_{B}) := \{u\in H^2(E_{\Gamma})\mid B(u)=0\}.$
\newline
{\rm{(ii)}} The operator $\bar{A}_B$ is positive definite and its spectrum
is bounded below by $\epsilon.$
\newline
{\rm{(iii)}} The operator
$$ (A_{\Gamma}, B) : C^{\infty}(E_{\Gamma}) \rightarrow
C^{\infty}(E_{\Gamma}) \oplus C^{\infty}(E_{\Gamma}
\mid _{\Gamma^{+} \sqcup \Gamma^{-}})
$$
defined by $(A_{\Gamma}, B)(u) = (A_{\Gamma}(u), B(u))$ can be
extended to an invertible operator $(A_{\Gamma}, B)\bar{ }$
$$(A_{\Gamma},B)\bar{ }: H^2(E_{\Gamma}) \rightarrow
L^2(E_{\Gamma}) \oplus H^{2 -\frac{1}{2}}(E_{\Gamma} \mid _{\Gamma^{+}
\sqcup \Gamma^{-}}).
$$
\endproclaim
\demo{Proof}
%We restrict ourselves to prove that $A_B$ is symmetric.
%Let $\{U_{\alpha}\}$ be an atlas of $M_{\Gamma}$ and $\{\rho_{\alpha}\}$ be a
%partition of unity subordinate to $\{U_{\alpha}\}$. For each $U_{\alpha}$, let
%$(x_{1},x_{2},\cdots,x_{d})$ be a local coordinate system and $(e_{1},e_{2},
%\cdots,e_{k})$ be an orthonormal frame of $\pi^{-1}(U_{\alpha})$. Then in each
%$U_{\alpha}$, $A$ can be expressed by $k\times k$ matrix of differential
%operators in $(x_{1},x_{2},\cdots,x_{k})$, say, $A=(a_{ij}(x))$, where
%$a_{ij}(x)$ is a differential operator of order 2.
%
%Let $u,v\in C^{\infty}(E_{\Gamma})$ with $B(u)=B(v)=0$.
%For each $U_{\alpha}$, choose a smooth function $\tau_{\alpha}$ such that
%$\tau_{\alpha}=1$ on $supp\rho_{\alpha}$ and $supp\tau_{\alpha}\subset
%U_{\alpha}$. Then
%\dddd
%$<A_{B}u,v>=\sum_{\alpha}<A_{B}\rho_{\alpha}u,v>=\sum_{\alpha}
%<A_{B}\rho_{\alpha}u,\tau_{\alpha}v>$  \newline
%\dddd
%$=\sum_{\alpha}\int_{U_{\alpha}}(\sum_{i,j}a_{ij}(x)\rho_{\alpha}u_{i}e_{i},
%\sum_{l}\tau_{\alpha}v_{l}e_{l})dvol$  \newline
%\dddd
%$=\sum_{\alpha}\int_{U_{\alpha}}\sum_{i,j}(a_{ij}(x)\rho_{\alpha}u_{i})
%(\tau_{\alpha}v_{i})dvol$  \newline
%\dddd
%$=\sum_{\alpha}\int_{U_{\alpha}}
%\sum_{i,j}(\rho_{\alpha}u_{i})(a_{ij}(x)^{\ast}
%\tau_{\alpha}v_{i})dvol$ \newline
%\dddd
%(by using integration by parts and the fact that $\rho_{\alpha}u_{i}$ and
%$\tau_{\alpha}u_{i}$ are 0 on the boundary of $U_{\alpha}$) \newline
%\dddd
%$=\sum_{\alpha}\int_{U_{\alpha}}(\rho_{\alpha}u,A^{\ast}(\tau_{\alpha}v))dvol
%=<u,A_{B}v>$,  since $A^{\ast}=A$. \newline
%\dddd
(i) Using a partition of unity and integration by parts one shows that $A_B$
is symmetric. Clearly, $\bar{A}_{B}$ is well defined and by a standard
argument self-adjoint.
To prove (ii) one first notices
that for any $u\in C^{\infty}(E_{\Gamma})$ with
$u\mid_{\Gamma^{+}\sqcup\Gamma^{-}}=0$,
one can find
a sequence $\{\phi_{n}\}$ such that supp$(\phi_{n})\subset M-\Gamma$ and
$\phi_{n}$ converges to $u$ in $H^{1}(E_{\Gamma})$.
Observe that $\langle A_{B}\phi_{n},\phi_{n}\rangle =
\langle A \phi_{n},\phi_{n} \rangle \geq \epsilon \| \phi_n\|^2$
and, integrating  by
parts, one concludes that
$$\langle A_{B}u,u \rangle =
\lim_{n\rightarrow\infty} \langle A\phi_{n},\phi_{n} \rangle \geq
\epsilon \|u\|^2.
$$
Thus (ii) follows.
\newline
(iii) As $A_{B}\bar{ }$ is injective,
so is the extension $(A_{\Gamma}, B)\bar{ }.$
To prove that this extension is onto, consider
$ f \in L^2(E_{\Gamma})$ and $\varphi \in
H^{2 - \frac{1}{2}}(E_{\Gamma}\mid _{\Gamma^{+}
\sqcup \Gamma^{-}}).$
Choose any section $v \in H^2(E_{\Gamma})$ so that $Bv =\varphi.$
As $\bar{A}_B$ is invertible, there exists $w \in
H^2(E_{\Gamma})$ satisfying $\bar{A}_Bw = f - \bar{A}_{\Gamma} v$
and the boundary conditions $Bw = 0.$ Therefore $u = w+v$ is an element
in $H^2(E_{\Gamma})$ with $(A_{\Gamma},B)\bar{ } u = (f, \varphi).$
Altogether one concludes that $(A_{\Gamma}, B)\bar{ }$ is
an isomorphism.
\qed
\enddemo

Set $\alpha_{k}=e^{i\frac{\pi+2k\pi}{d}}$ for $0\leq k\leq d-1,$
where $d=dim(M).$

\proclaim{Lemma 3.2}
The following operators are invertible for
$0\leq k\leq d-1$ and $t \geq 0$
$$(A_{\Gamma}-\alpha_{k}t,B):C^{\infty}(E_{\Gamma})\rightarrow
C^{\infty}(E_{\Gamma})\oplus
C^{\infty}(E_{\Gamma}\mid_{\Gamma^{+}\sqcup\Gamma^{-}}).$$
\endproclaim
\demo{Proof}
As $\alpha_{k} \in \Bbb C \setminus \Bbb R ^{+}$ and thus, for
$t \geq 0,$ $\alpha_{k}t\not\in Spec(A_{B}),$
the operator $(A_{\Gamma}-\alpha_{k}t,B)$ is
injective.
To prove that this operator is onto one argues as in the proof
of Lemma 3.1 (iii).
\qed
\enddemo

Since $(A_{\Gamma}-\alpha_{k}t,B)$ is invertible, we can define the Poisson
operator $P(\alpha_{k} t)$ associated to $(A_{\Gamma}-\alpha_{k}t,B),$
$P(\alpha_k t): C^{\infty}(E_{\Gamma} \mid_{\Gamma^{+} \sqcup \Gamma^{-}})
\rightarrow C^{\infty}(E_{\Gamma}),$
i.e. for $\varphi \in C^{\infty}(E_{\Gamma}
\mid_{\Gamma^{+} \sqcup \Gamma^{-}}),$
$u = P(\alpha_k t) \varphi$
is the solution in $C^{\infty}(E_{\Gamma})$ of
$(A_{\Gamma} - \alpha_{k} t) u =0$ with boundary conditions
$u\mid_{\Gamma^{+} \sqcup \Gamma^{-}}
= \varphi.$
\newline
Let $R(\alpha_k t) : C^{\infty}(E\mid _\Gamma) \rightarrow
C^{\infty}(E \mid_{\Gamma})$ be the Dirichlet to Neumann operator
corresponding to $A_{\Gamma} - \alpha_k t , B$ and $C.$ Then the following
result holds:
\proclaim{Lemma 3.3}
For $0\leq k\leq d-1$, and $t\geq 0$, $R(\alpha_{k}t)$ is an invertible
classical $\Psi$DO of order 1, which is elliptic with parameter $t$ of
weight 1.
\endproclaim
\demo{Proof}
In a sufficiently small collar neighborhood $U$ of $\Gamma,$ choose
coordinates $x=(x',s)$ such that
$(x',0)\in\Gamma$ and
$\frac{\partial}{\partial s}\mid_{(x',0)}= \nu_{(x',0)}$. Let $\xi =
(\xi',\eta)$ be coordinates in the cotangent space corresponding to
the coordinates $(x',s).$
Let $D_{s}=\frac{1}{i}\frac{\partial}{\partial s}$ and write
$(A-\alpha_{k}t)=A_{2}D_{s}^{2}+A_{1}D_{s}+A_{0}$, where the
$A_{j}$'s are differential operators
of order at most $2-j$. The $A_{j}$'s induce, when restricted to
$\Gamma$, differential operators, again denoted by $A_{j},$
$A_{j} :C^{\infty}(E\mid_{\Gamma})
\rightarrow C^{\infty}(E\mid_{\Gamma})$.
Since $\sigma_{L}(x,(\xi',\eta))=\parallel(\xi',\eta)\parallel^{2}$
 and since $\nu_{(x',0)}$ is the unit normal to $\Gamma$ at
$(x',0),$ one has $A_{2}(x)=Id_{x} \in End_{x}(E_x,E_x)$ on $\Gamma$.

For any $\varphi\in C^{\infty}(E\mid_{\Gamma})$ and $t \geq 0,$
 we can choose $u\in C^{\infty}
(E_{\Gamma})\cap C(E)$ such that $(A-\alpha_{k}t)u=0$ on $M-\Gamma$ and
$u\mid_{\Gamma^{+}}=\varphi= u\mid_{\Gamma^{-}}$.
% Note that $(A-\alpha_{k}t)=
%-A_{2}\frac{d^{2}}{ds^{2}}+\frac{1}{i}A_{1}\frac{d}{ds}
%+A_{0}$.
Then $\frac{\partial u}{\partial s}(x',s)$ has a jump across $\Gamma,$
which is
$-R(\alpha_{k}t)(\varphi)(x')$.
Hence
$$\frac{\partial u}{\partial s}(x',s)
= - R(\alpha_{k}t)(\varphi)(x') H(s)+v(x',s),
$$
 where $v(x',s)\in C^{\infty}(E_{\Gamma}\mid _{U})\cap C(E \mid _{U})$
 and $H(s)$ is the
Heavyside function. Therefore, on $U,$
$$(A-\alpha_{k}t)u=A_{2}R(\alpha_{k}t)(\varphi)\otimes\delta_{\Gamma}
-A_{2}\frac{\partial v}{\partial s}+\frac{1}{i}A_{1}
\frac{\partial u}{\partial s}+A_{0}u.$$
%So, again on $U,$
%$$(A-\alpha_{k}t)u=A_{2}\cdot(\cdot\otimes\delta_{\Gamma})\cdot
%R(\alpha_{k}t)(\varphi)+(-A_{2}\frac{\partial v}{\partial s}+
%\frac{1}{i}A_{1}\frac{\partial u}{\partial s}+A_{0}).
%$$
Since $(A-\alpha_{k}t)u=0$ on $M-\Gamma$, we conclude that, on
$U \cap (M \setminus \Gamma),$
$$-A_{2}\frac{\partial v}{\partial s}+
\frac{1}{i}A_{1}\frac{\partial u}{\partial s}+A_{0}u
=0.$$
As $-A_{2}\frac{\partial v}{\partial s}+\frac{1}{i}A_{1}
\frac{\partial u}{\partial s}+A_{0}u
\in L^{2}(E \mid _{U})$, it follows that
$$(A-\alpha_{k}t)u=A_{2}\cdot(\cdot\otimes\delta_{\Gamma})\cdot R(\alpha_{k}t)
\varphi.
$$
%Hence
%$Id=J\cdot (A-\alpha_{k}t)^{-1}\cdot A_{2}\cdot
%(\cdot\otimes\delta_{\Gamma})\cdot R(\alpha_{k}t)$,
Using that
 $A_{2}=Id$ on $\Gamma$,
one therefore obtains
$Id=J\cdot (A-\alpha_{k}t)^{-1}\cdot(\cdot\otimes\delta_{\Gamma})\cdot
R(\alpha_{k}t)$
where
$J$ is the restriction operator to $\Gamma$.
 From this identity
it follows that $R(\alpha_{k}t)$ is invertible. Moreover, setting $\phi=
R(\alpha_{k}t)\varphi$,
$$R(\alpha_{k}t)^{-1}\phi=J\cdot (A-\alpha_{k}t)^{-1}\cdot
(\phi\otimes\delta_{\Gamma})$$
$$=J\cdot\int_{\Bbb R^{d-1}}\int_{\Bbb R}e^{i(x',s)\cdot (\xi',\eta)}
[(A-\alpha_{k}t)^{-1}
(\phi\otimes\delta_{\Gamma})]\hat{\, }(\xi',\eta)d\eta d\xi'$$
$$=\int_{\Bbb R^{d-1}}
e^{i x'\cdot\xi'}\int_{\Bbb R}\sigma((A-\alpha_{k}t)^{-1})
(x',0,\xi',\eta)\hat{\phi}(\xi')\cdot\frac{1}{\sqrt{2\pi}}d\eta d\xi'.$$
Hence $R(\alpha_{k}t)^{-1}$ is a classical $\Psi$DO of order -1 with
symbol
$$\frac{1}{\sqrt{2\pi}}\int_{\Bbb R}\sigma((A-\alpha_{k}t)^{-1})
(x',0,\xi',\eta)d\eta,
$$
and therefore $R(\alpha_{k}t)$ is a classical $\Psi$DO of order 1 with
parameter t of weight 1.
The ellipticity with parameter of $R(\alpha_kt)$
follows from the explicit formula of the symbol.
\qed
\enddemo

\proclaim{Lemma 3.4}
For $\epsilon' <\frac{\pi}{d}$ sufficiently small and $0 \leq k \leq d-1,$
the operator $R(\alpha_{k}t)$ does not have any eigenvalues
in $\Lambda_{\epsilon'}$, where $\Lambda_{\epsilon'}=
\{z\in\Bbb C\mid \pi-\epsilon' < arg(z)<\pi+\epsilon' \text{ or }|z| <
\epsilon' \}.$ Hence
$R(\alpha_{k}t)$ has $\pi$ as an Agmon angle.
\endproclaim
\demo{Proof}
By assumption, $A:C^{\infty}(E)
\rightarrow C^{\infty}(E)$ is essentially self-adjoint and positive definite.
Let $\{\psi_{j}\}_{j\geq1}$ be a complete orthonormal system of
eigensections of $A$
with corresponding eigenvalues $\{\lambda_{j}\}_{j\geq 1}.$
Then $(A-\alpha_{k}t)^{-1}
\psi_{j}=(\lambda_{j}-\alpha_{k}t)^{-1}\psi_{j}$.

Moreover, for any
$\varphi\in C^{\infty}(E\mid_{\Gamma}),\varphi\otimes\delta_{\Gamma}$
is an element in
$H^{-1}(E)$ and $(A-\alpha_{k}t)^{-1}
(\varphi\otimes\delta_{\Gamma})\in L^{2}(E).$
Therefore
$$\langle (A-\alpha_{k}t)^{-1}\varphi\otimes\delta_{\Gamma},\psi_{j} \rangle
=$$
$$\langle \varphi\otimes
\delta_{\Gamma},((A-\alpha_{k}t)^{-1})^{\ast}\psi_{j} \rangle
%=(\lambda_{j}-
%\alpha_{k}t)^{-1}$
%$\langle \varphi\otimes\delta_{\Gamma},
%\psi_{j}\rangle
=(\lambda_{j}-\alpha_{k}t)^{-1}\int_
{\Gamma}(\varphi,\psi_{j})d\mu_{\Gamma},$$ where $d\mu_{\Gamma}$ is the volume
  form on $\Gamma$ induced from the metric on $M$ and $(\cdot,\cdot)$ is the
Hermitian inner product on $E$.

Since
$R(\alpha _{k}t)^{-1}=J\cdot (A-\alpha_{k}t)^{-1}\cdot (\cdot\otimes\delta_
{\Gamma}),$ one obtains for
$\varphi_{1}, \varphi_{2} \in C^{\infty}(E\mid _{\Gamma})$
$$\langle R(\alpha_{k}t)^{-1}\varphi_{1},\varphi_{2} \rangle
=\sum^{\infty}_{j=1}(\lambda_{j}-\alpha_{k}t)^{-1}
\int_{\Gamma}(\varphi_{1},\psi_{j})d\mu_{\Gamma}
\int_{\Gamma}(\psi_{j}, \varphi_{2})d\mu_{\Gamma}.$$
Together with Lemma 3.3 this implies that $\Lambda_{\epsilon'}$
has an empty intersection with $SpecR(\alpha_{k} t).$
\qed
\enddemo

Using the above formula, one obtains as an immediate consequence the following
\proclaim{Corollary 3.5}
The operator $R = R(0)$ is essentially self-adjoint and
positive definite.
\endproclaim
%\demo{Proof}
%The result follows from the fact that
%$$<R^{-1}\varphi_{1},\varphi_{2}>=\sum_{j=1}^{\infty}\lambda_{j}^{-1}
%\int_{\Gamma}(\varphi_{1},\psi_{j})d\mu_{\Gamma}\int_{\Gamma}(\psi_{j},
%\varphi_{2})d\mu_{\Gamma}.$$
%\qed
%\enddemo

Next we are collecting a number of results about operators involving the
$d-$th power of $A$ and the submanifold $\Gamma$.

Consider the families of operators $A^{d}+t^{d}$ and
$A_{\Gamma}^{d}+t^{d}$ for nonnegative real numbers $t$.
Then $A^{d}+t^{d}$ and
$A_{\Gamma}^{d}+t^{d}$ are elliptic differential operators with parameter,
where the weight of t is 2. Note that
$$A_{\Gamma}^{d}+t^{d}=(A_{\Gamma}-te^{i\frac{\pi}{d}})(A_{\Gamma}-
te^{i\frac{3\pi}{d}})\cdots (A_{\Gamma}-te^{i\frac{\pi+2\pi(d-1)}{d}}).$$

Let us introduce the boundary conditions $B_{d}(t), C_{d}(t)$ by setting
$$B_{d}(t)=(B,B(A_{\Gamma}-\alpha_{0}t),B(A_{\Gamma}-\alpha_{1}t)(A_{\Gamma}-
\alpha_{0}t), \cdots ,
B(A_{\Gamma}-\alpha_{d-2}t)\cdots (A_{\Gamma}-\alpha_{0}t)),$$
and
$$C_{d}(t)=(C,C(A_{\Gamma}-\alpha_{0}t),C(A_{\Gamma}-\alpha_{1}t)(A_{\Gamma}-
\alpha_{0}t), \cdots ,
C(A_{\Gamma}-\alpha_{d-2}t)\cdots (A_{\Gamma}-\alpha_{0}t)).$$
It follows from Lemma 3.2 that the following operator is invertible
$$(A_{\Gamma}^{d}+t^{d},B_{d}(t)):C^{\infty}(E_{\Gamma})\rightarrow
C^{\infty}(E_{\Gamma})\oplus \left(\oplus_{d}
C^{\infty}(E_{\Gamma}\mid_{\Gamma^{+}\sqcup\Gamma^{-}})\right).$$
Therefore the corresponding  Poisson operator
$\tilde{P_{d}}(t) : \oplus _{d} C^{\infty}(E\mid _{\Gamma^{+}
\sqcup \Gamma^{-}})
\rightarrow
C^{\infty}(E_{\Gamma})$
is well defined.

\proclaim{Lemma 3.6}
The Poisson operator $\tilde{P_{d}}(t)$ associated to
$(A_{\Gamma}^{d}+t^{d},B_{d}(t))$ is given by
$$\tilde{P_{d}}(t)(\varphi_{0},\cdots,\varphi_{d-1})
=P(\alpha_{0}t)\varphi_{0}+ (A_{\Gamma}-\alpha_{0}t)_{B}^{-1}P(\alpha_{1}t)
\varphi_{1}+\cdots +$$
$$(A_{\Gamma}-\alpha_{0}t)_{B}^{-1}(A_{\Gamma}-\alpha_{1}t)_{B}
^{-1}\cdots (A_{\Gamma}-\alpha_{d-2}t)_{B}^{-1}P(\alpha_{d-1}t)\varphi_{d-1},$$
where $(A_{\Gamma}-\alpha_{k}t)_{B}$ is the restriction of
$A_{\Gamma}-\alpha_{k}t$ to $\{u\in C^{\infty}(E_{\Gamma})\mid Bu=0 \}$.
\endproclaim
\demo{Proof}
Denoting the right hand side of the claimed identity by
$Q_{d}(t)(\varphi_{0},\cdots,\varphi_{d-1})$ one obtains
$$(A_{\Gamma}^{d}+t^{d})\cdot Q_{d} (t)
(\varphi_{0},\cdots,\varphi_{d-1})=0.$$
Moreover, for $0 \leq k \leq d-1, Q_{d}(t)(\varphi_{0},\cdots,\varphi_{d-1})$
satisfies the boundary conditions
$$(B(A_{\Gamma}-\alpha_{k-1}t)(A_{\Gamma}-\alpha_{k-2}t)\cdots
(A_{\Gamma}-
\alpha_{0}t)) Q_{d}(t)(\varphi_{0},\cdots,\varphi_{d-1})=
$$
$$B(A_{\Gamma}-\alpha_{k-1}t)\cdots (A_{\Gamma}-\alpha_{0}t)P(\alpha_{0}t)
\varphi_{0}+\cdots +
$$
$$B(A_{\Gamma}-\alpha_{k-1}t)\cdots (A_{\Gamma}-\alpha_{0}t)
(A_{\Gamma}-\alpha_{0}t)_{B}^{-1}\cdots
(A_{\Gamma}-\alpha_{k-1}t)_{B}^{-1}P(\alpha_{k}t)\varphi_{k}+\cdots +
$$
$$B(A_{\Gamma}-\alpha_{k-1}t)\cdots
(A_{\Gamma}-\alpha_{0}t)(A_{\Gamma}-\alpha_{0}t)_{B}^{-1}\cdots
(A_{\Gamma}-\alpha_{d-2}t)_{B}^{-1}P(\alpha_{d-1}t)\varphi_{d-1}
$$
$$=\varphi_{k},
$$
since $(A_{\Gamma}-\alpha_{j}t)P(\alpha_{j}t)=0$ and
$B(A_{\Gamma}-\alpha_{j}t)_{B}^{-1}=0.$ These two properties of $Q_{d}(t)$
establish the claimed identity.
\qed
\enddemo

Further let us consider the boundary conditions $B_{d}(t)$ and $C_{d}(t)$ for
$t=0$. Note that
$$B_{d}(0)=(B,BA_{\Gamma},\cdots,BA_{\Gamma}^{d-1});
C_{d}(0)=(C,CA_{\Gamma},\cdots,CA_{\Gamma}^{d-1}).$$
Let $\Omega(t)$ be the following lowertriangular  $d\times d$ matrix
$$
\Omega(t)=\left(\matrix
1 & 0 & \cdots & 0 \\
\alpha_{0}t & 1 & \cdots & 0 \\
\vdots & \vdots & \ddots & \vdots \\
\alpha_{0}^{d-1}t^{d-1} & t^{d-2}\sum^{d-2}_{k=0}\alpha_{0}^{d-2-k}
\alpha_{1}^{k} & \cdots & 1
\endmatrix \right).
$$
Then $B_{d}(0)=\Omega(t)B_{d}(t)$ as well as $C_{d}(0)=\Omega(t)C_{d}(t).$
Let $P_{d}(t):=\tilde{P}_{d}(t)\Omega(t)^{-1}$  and notice that $P_{d}(t)$ is
the Poisson operator corresponding to $(A_{\Gamma}^{d}+t^{d},B_{d}(0)).$

Consider the Dirichlet to Neumann operator $\tilde{R_{d}}(t)=
\bigtriangleup_{if}\cdot C_{d}(t)
\cdot \tilde{P_{d}}(t)\cdot \bigtriangleup_{ia}$
corresponding to
$A_{\Gamma}^{d}+t^{d}, B_{d}(t)$ and $C_{d}(t).$ Then
$$\tilde{R_{d}}(t)(\varphi_{0},\cdots,\varphi_{d-1})=$$

$\bigtriangleup_{if}\cdot
\left(C,C(A_{\Gamma}-\alpha_{0}t),\cdots,C(A_{\Gamma}-\alpha_{d-2}t)\cdots
(A_{\Gamma}-\alpha_{0}t)\right)\cdot$ \newline
\dddd
$(P(\alpha_{0}t)
\bigtriangleup_{ia}\varphi_{0}+(A_{\Gamma}-\alpha_{0}t)_{B}^{-1}
P(\alpha_{1}t)\bigtriangleup_{ia}
\varphi_{1}+ \cdots +$ \newline
\dddd
$(A_{\Gamma}-\alpha_{0}t)_{B}^{-1}(A_{\Gamma}-\alpha_{1}t)_{B}^{-1}\cdots
(A_{\Gamma}-\alpha_{d-2}t)_{B}^{-1}P(\alpha_{d-1}t)
\bigtriangleup_{ia}\varphi_{d-1}).$ \newline
\dddd

Thus $\tilde{R_{d}}(t):
\oplus_{d}C^{\infty}(E\mid_{\Gamma})\rightarrow \oplus_{d}
C^{\infty}(E\mid_{\Gamma})$ can be represented by a $d\times d$ matrix of
upper triangular form,
$$
\left(\smallmatrix
R(\alpha_{0}t) &\bigtriangleup_{if}C(A_{\Gamma}-\alpha_{0}t)_{B}^{-1}
P(\alpha_{1}t)\bigtriangleup_{ia} &\cdots &\bigtriangleup_{if}
C(A_{\Gamma}-\alpha_{0}t)_{B}^{-1}\cdots
(A_{\Gamma}-\alpha_{d-2}t)_{B}^{-1}P(\alpha_{d-1}t)
\bigtriangleup_{ia}\\
0 & R(\alpha_{1}t) &\cdots &\bigtriangleup_{if}
C(A_{\Gamma}-\alpha_{1}t)_{B}^{-1}\cdots (A_{\Gamma}-\alpha_{d-2}t)_{B}^{-1}
P(\alpha_{d-1}t)\bigtriangleup_{ia}\\
\vdots &\vdots &\ddots &\vdots\\
0 & 0 &\cdots & R(\alpha_{d-1}t)
\endsmallmatrix \right),
$$
where $R(\alpha_{k}t)$ is the Dirichlet to Neumann operator corresponding to
$A_{\Gamma}-\alpha_{k}t , B$ and $C$ defined earlier. In particular,
we conclude that $\tilde{R_{d}}(t)$ is invertible and has $\pi$ as an Agmon
angle.

Finally introduce
the Dirichlet to Neumann operator
$R_{d}(t)$
associated to
$A_{\Gamma}^{d}+t^{d}, B_{d}(0)$ and $C_{d}(0).$ Then
\newline
$$\tilde{R_{d}}(t)=\bigtriangleup_{if}\cdot C_{d}(t)\cdot \tilde{P_{d}}(t)\cdot
\bigtriangleup_{ia}
=\bigtriangleup_{if}\cdot \Omega(t)^{-1}\cdot C_{d}(0)\cdot P_{d}(t)\cdot
\Omega(t)\cdot \bigtriangleup_{ia}$$
$$=\Omega(t)^{-1}\cdot \bigtriangleup_{if}\cdot C_{d}(0)\cdot P_{d}(t)\cdot
\bigtriangleup_{ia}\cdot \Omega(t)
=\Omega(t)^{-1}\cdot R_{d}(t)\cdot \Omega(t).$$
As a consequence, $R_{d}(t)$ has the same spectrum as $\tilde{R_{d}}(t)$ and
therefore, $R_{d}(t)$ is invertible, has $\pi$ as an Agmon angle and
satisfies $\log Det(R_{d}(t))= \log Det (\tilde{R_{d}}(t)).$
In view of the fact that
$\tilde{R_{d}}(t)$ is of  upper triangular form one has
$$logDet(\tilde{R_{d}}(t))=\sum^{d-1}_{k=0}logDet(R(\alpha_{k}t))
%=logDet(R_{d}(t))
.$$

As $A$ is positive and essentially selfadjoint, the operator
$A^{d}+t^{d}:C^{\infty}(E)\rightarrow C^{\infty}(E)$ is
invertible for $t \geq 0.$ Using the kernel
$k_{t}(x,y)$ of $(A^{d}+t^{d})^{-1}$ this operator
can be extended to
$C^{\infty}(E_{\Gamma})$ by setting ($u \in C^{\infty}(E_{\Gamma})$)
 $${((A^{d}+t^{d})^{-1})}_{\Gamma}u(x)=\int_{M_{\Gamma}}k_{t}(x,y)u(y)dy.$$
It follows from Lemma 3.2 that
$$(A^d_{\Gamma}+t^d,B_d(t)):C^{\infty}(E_{\Gamma})\rightarrow
C^{\infty}(E_{\Gamma})\oplus \left(\oplus_{d}
C^{\infty}(E_{\Gamma}\mid_{\Gamma^{+}\sqcup\Gamma^{-}}) \right).$$
is invertible. Thus, since
$B_{d}(0)=\Omega(t)B_{d}(t),$ we conclude that $(A_{\Gamma}^d+t^d,B_d(0))$
is invertible as well.
Denote by $(A_{\Gamma}^{d}+t^{d})_{B_{d}(0)}$ the restriction of
$A^{d}_{\Gamma}+t^{d}$ to \newline  $\{u\in C^{\infty}
(E_{\Gamma})\mid B_{d}(0)u=0\}$
and let $(A^{d}_{\Gamma}+t^{d})_{B_{d}(0)}^{-1}$
be its inverse.

\proclaim{Lemma 3.7}
$(A_{\Gamma}^{d}+t^{d})_{B_{d}(0)}^{-1}=((A^{d}+t^{d})^{-1})_{\Gamma}-
P_{d}(t)\cdot B_{d}(0)\cdot ((A^{d}+t^{d})^{-1})_{\Gamma}$
\endproclaim
\demo{Proof}
%Since $({A^{d}}_{\Gamma}+t^{d})\cdot {(A^{d}+t^{d})^{-1}}_{\Gamma}=Id$ ,
Denote by $Q(t)$ the right hand side of the claimed identity. One verifies
that for $u \in C^{\infty}(E_{\Gamma})$
$$
(A^{d}_{\Gamma}+t^{d})Q(t)u = u
$$
and
$$
B_{d}(0)Q(t)u =
B_{d}(0)\cdot ((A^{d} +t^{d})^{-1})_{\Gamma} u -
B_{d}(0)\cdot ((A^{d} +t^{d})^{-1})_{\Gamma} u =0.
$$
These two identities imply that
$Q(t) = (A^d_{\Gamma} + t^{d} )^{-1}_{B_{d}(0)}.$
\qed
\enddemo

\proclaim{Lemma 3.8}
{\rm{(i)}} $\frac{d}{dt}P_{d}(t)=-dt^{d-1}
(A^{d}_{\Gamma}+t^{d})_{B_{d}(0)}^{-1}\cdot
P_{d}(t)$ \newline
{\rm{(ii)}} $R_{d}(t)^{-1}\cdot\frac{d}{dt}R_{d}(t)=$
$-dt^{d-1}R_{d}(t)^{-1}\cdot\bigtriangleup_{if}
\cdot C_{d}(0)\cdot (A_{\Gamma}^{d}+t^{d})_{B_{d}(0)}^{-1}\cdot P_{d}(t)\cdot
\bigtriangleup_{ia}.$  \newline
In particular, d being the dimension of $M,$
$R_{d}(t)^{-1}\cdot\frac{d}{dt}R_{d}(t)$ is of trace class.
\endproclaim
\demo{Proof}
(i) Derive $(A_{\Gamma}^{d}+t^{d})\cdot P_{d}(t)=0$ with respect to t to
obtain
$$(A_{\Gamma}^{d}+t^{d})\cdot\frac{d}{dt}
P_{d}(t)=-\frac{d}{dt}(A_{\Gamma}^{d}+t^{d})\cdot P_{d}(t)=-dt^{d-1}
P_{d}(t).
$$
Similarly, deriving $B_{d}(0)\cdot P_{d}(t)=Id$ with respect to t yields
$B_{d}(0)\frac{d}{dt}P_{d}(t)=0.$
Hence
$$(A^{d}_{\Gamma}+t^{d})_{B_{d}(0)}\cdot \frac{d}{dt}P_{d}(t)=-dt^{d-1}
P_{d}(t)
$$
and therefore
$$\frac{d}{dt}P_{d}(t)=-dt^{d-1}(A^{d}_{\Gamma}+t^{d})_{B_{d}(0)}^{-1}\cdot
P_{d}(t).
$$
(ii) follows from the definition of $R_{d}(t)$ and (i).
\qed
\enddemo

Taking into account that
$\bigtriangleup_{ia}(\oplus_{d}C^{\infty}(E\mid_{\Gamma}))=
\{(\varphi,\varphi)\mid \varphi\in\oplus_{d}C^{\infty}(E\mid_{\Gamma})\}$
we may define
$Pr_{\Gamma}:\bigtriangleup_{ia}(\oplus_{d}C^{\infty}
(E\mid_{\Gamma}))\rightarrow \oplus_{d}C^{\infty}(E\mid_{\Gamma})$ by
$Pr_{\Gamma}(\varphi,\varphi)=\varphi$.

\proclaim{Corollary 3.9}
$$R_{d}(t)^{-1}\cdot\frac{d}{dt}R_{d}(t)=
dt^{d-1}Pr_{\Gamma}\cdot B_{d}(0)\cdot
((A^{d}+t^{d})^{-1})_{\Gamma}\cdot P_{d}(t)\cdot\bigtriangleup_{ia}$$
\endproclaim
\demo{Proof}
By Lemma 3.7 and 3.8
$$R_{d}(t)^{-1}\cdot\frac{d}{dt}R_{d}(t)
=-dt^{d-1}(\bigtriangleup_{if}\cdot C_{d}(0)\cdot P_{d}(t)\cdot
\bigtriangleup_{ia})^{-1}\cdot
$$
$$\bigtriangleup_{if}\cdot C_{d}(0)
\cdot \left( ((A^{d}+t^{d})^{-1})_{\Gamma}-P_{d}(t)\cdot B_{d}(0)\cdot
((A^{d}+t^{d})^{-1})_{\Gamma} \right) \cdot P_{d}(t)\cdot\bigtriangleup_{ia}.
$$
%Next consider $\bigtriangleup_{if}\cdot C_{d}(0)\cdot
%((A^{d}+t^{d})^{-1})_{\Gamma}\cdot
%P_{d}(t)\cdot\bigtriangleup_{ia}:\oplus_{d}C^{\infty}(E\mid_{\Gamma})
%\rightarrow\oplus_{d}C^{\infty}(E\mid_{\Gamma})$. \newline
%For any $\varphi\in \oplus_{d}C^{\infty}(E\mid_{\Gamma}), P_{d}(t)\cdot
%\bigtriangleup_{ia}(\varphi)\in (C^{\infty}(E_{\Gamma})\cap C^{0}(E))
%\subset L^{2}(E).$
%Hence\newline
%$(\bigtriangleup_{if}\cdot C_{d}(0)\cdot
%(A^{d}+t^{d})^{-1}_{\Gamma})\mid_{H^{0}(E)}$ is the composition of
%the following maps
%$$ H^{0}(E) @>(A^{d}+t^{d})^{-1}_{\Gamma}>> H^{2d}(E)@>C_{d}(0)>>
%H^{2d-1-\frac{1}{2}}(E\mid_{\Gamma})\oplus H^{2d-3-\frac{1}{2}}
%(E\mid_{\Gamma})\oplus \cdots
%\oplus H^{\frac{1}{2}}(E\mid_{\Gamma})$$
%$$@>\bigtriangleup_{if} >>
%H^{2d-\frac{3}{2}}(E\mid_{\Gamma})\oplus \cdots\oplus H^{\frac{1}{2}}
%(E\mid_{\Gamma}).$$
Clearly $\bigtriangleup_{if}\cdot C_{d}(0)\cdot ((A^{d}+t^{d})^{-1})_{\Gamma}
(u)=0$ for $u\in C^{\infty}(E)$
%Since $C^{\infty}(E)$ is dense in $H^{0}(E), \bigtriangleup_{if}\cdot C_{d}(0)
%\cdot (A^{d}+t^{d})^{-1}_{\Gamma}=0$ on $H^{0}(E)$ \newline
\newline
and thus
$\bigtriangleup_{if}
\cdot C_{d}(0)\cdot ((A^{d}+t^{d})^{-1})_{\Gamma}\cdot P_{d}(t)\cdot
\bigtriangleup_{ia}=0.
$
Therefore $R_{d}(t)^{-1}\cdot\frac{d}{dt}R_{d}(t)=
dt^{d-1}(\bigtriangleup_{if}\cdot C_{d}(0)\cdot
P_{d}(t)\cdot\bigtriangleup_{ia})^{-1}\cdot\bigtriangleup_{if}\cdot
C_{d}(0)\cdot P_{d}(t)\cdot B_{d}(0)\cdot ((A^{d}+t^{d})^{-1})_{\Gamma}\cdot
P_{d}(t)\cdot\bigtriangleup_{ia}.$
\newline
Note that for any $u\in C^{\infty}(E_{\Gamma}),$
the boundary values of $((A^{d}+t^{d})^{-1})_{\Gamma}u$
on $\Gamma^{+}$ and $\Gamma^{-}$ are the same,
i.e.
$$
B_{d}(0)((A^{d}+t^{d})^{-1})_{\Gamma}u\mid_{\Gamma^{+}}=
B_{d}(0)((A^{d}+t^{d})^{-1})_{\Gamma}u\mid_{\Gamma^{-}}.
$$
Hence
$$B_{d}(0)\cdot ((A^{d}+t^{d})^{-1})_{\Gamma}\cdot P_{d}(t)\cdot
\bigtriangleup_{ia}=\bigtriangleup_{ia}\cdot Pr_{\Gamma}\cdot B_{d}(0)\cdot
((A^{d}+t^{d})^{-1})_{\Gamma}\cdot P_{d}(t)\cdot\bigtriangleup_{ia}.
$$
\qed
\enddemo

As $(A^{d}+t^{d})^{-1}, (A^{d}_{\Gamma}+t^{d})^{-1}_{B_{d}(0)}$ and
$R_{d}(t)^{-1}\frac{d}{dt}R_{d}(t)$ are of trace class
we can apply the well known variational formula for regularized determinants:

\proclaim{Lemma 3.10}
Let $Q(t)$ denote any of the operators $A^{d}+t^{d}$,
$(A^{d}+t^{d})_{B_{d}(0)}$ or $R_{d}(t)$.
Then, for any $t \geq 0,$
$$\frac{d}{dt}
logDetQ(t)=tr(Q(t)^{-1}\frac{d}{dt}Q(t)).$$
\endproclaim

\dd
{\bf 4. Proof of the Theorem 1.1}
\ddd
In the case where the operator $A^{-1}$ is of trace class, the proof of
Theorem 1.1 is considerably simpler. Unfortunately, this is only the case
if the dimension d of $M$ is equal to 1.
Our strategy is to first prove a version of Theorem 1.1 for $A^{d}$
(Lemma 4.1), using the fact that $(A^{d})^{-1}$ is of trace class.
Together with the auxiliary results of section 3
and the asymptotic expansion derived in section 2, the proof of Theorem 1.1
is then completed.

\proclaim{Lemma 4.1}
Let $A^{d}+t^{d}$ and $(A_{\Gamma}^{d}+t^{d},B_{d}(0))$ be as above.
Then, for $t \geq 0,$
$$\frac{d}{dt}\left(logDet(A^{d}+t^{d})-
logDet( A_{\Gamma}^{d}+t^{d},B_{d}(0))\right)
= \frac{d}{dt}logDetR_{d}(t).$$
\endproclaim
\demo{Proof}
Define $w(t):=\frac{d}{dt}(logDet(A^{d}+t^{d})-
logDet((A_{\Gamma}^{d}+t^{d}),B_{d}(0)))$. By Lemma 3.10 and Lemma 3.7
\newline
\dddd
$w(t)=tr(\frac{d}{dt}(A^{d}+t^{d})\cdot (A^{d}+t^{d})^{-1}
-\frac{d}{dt}(A_{\Gamma}^d+t^{d})_{B_{d}(0)}\cdot
(A_{\Gamma}^{d}+t^{d})_{B_{d}(0)}^{-1})$ \newline
\dddd
$=dt^{d-1}tr((A^{d}+t^{d})^{-1}-
(A_{\Gamma}^{d}+t^{d})_{B_{d}(0)}^{-1})$
\newline
\dddd
$=dt^{d-1}tr(((A^{d}+t^{d})^{-1})_{\Gamma}-
(A_{\Gamma}^{d}+t^{d})^{-1}_{B_{d}(0)})$
\newline
\dddd
$=dt^{d-1}tr(P_{d}(t)
\cdot B_{d}(0)\cdot ((A^{d}+t^{d})^{-1})_{\Gamma}).$
\newline
\dddd

On the other hand, by Lemma 3.10, Corollary 3.9 and the commutativity
of the trace,
$$\frac{d}{dt}logDetR_{d}(t)=tr(\frac{d}{dt}R_{d}(t)\cdot R_{d}(t)
^{-1})$$
$$=dt^{d-1}tr(Pr_{\Gamma}\cdot B_{d}(0)\cdot ((A^{d}+t^{d})^{-1})_{\Gamma}
\cdot  P_{d}(t)\cdot\bigtriangleup_{ia})$$
$$=dt^{d-1}tr(P_{d}(t)\cdot
\bigtriangleup_{ia}\cdot Pr_{\Gamma}\cdot B_{d}(0)\cdot
((A^{d}+t^{d})^{-1})_{\Gamma})$$
$$=dt^{d-1}tr(P_{d}(t)\cdot B_{d}(0)\cdot ((A^{d}+t^{d})^{-1})_{\Gamma}).$$
Combining the above two identities shows that
$$ w(t) = \frac{d}{dt} log DetR_{d}(t).$$
\qed
\enddemo

Since $logDetR_{d}(t)=\sum^{d-1}_{k=0}logDetR(\alpha_{k}t)$, we
conclude from Lemma 4.1 that
$$logDet(A^{d}+t^{d})-logDet((A_{\Gamma}^{d}+t^{d}),B_{d}(0))= \tilde{c}
+\sum^{d-1}_{k=0}logDetR(\alpha_{k}t), \eqno (4.1)$$
where $\tilde{c}$ is independent of $t$.

Note that $logDet(A^{d}+t^{d})$, $logDet(A_{\Gamma}^{d}+t^{d},B_{d}(0))$ and
$logDetR(\alpha_{k}t)
(0\leq k\leq d-1)$ have asymptotic expansions as $t\rightarrow +\infty$.
Since the eigenvalues of $A^{d}+t^{d}$ and $(A_{\Gamma}^{d}+t^{d})_{B_{d}(0)}$
satisfy
the condition in Proposition 2.7, the constant terms in the asymptotic
expansions of
$logDet(A^{d}+t^{d})$ and $logDet((A_{\Gamma}^{d}+t^{d}),B_{d}(0))$ are zero.
Let $\pi_{0}(R(\alpha_{k}t))$ be the constant term in the asymptotic
expansion of
$logDet(R(\alpha_{k}t))$. Then $\tilde{c}=-\sum^{d-1}_{k=0}
\pi_{0}(R(\alpha_{k}t)),$
which is computable in terms of the symbol of
$R(\alpha_{k}t)$ by Theorem 2.4.

\proclaim{Lemma 4.2}
{\rm{(i)}} $Det(A_{\Gamma}^{d},B_{d}(0))=(Det(A_{\Gamma},B))^{d};$
{\rm{(ii)}} $Det(A^{d})=(DetA)^{d}.$
\endproclaim
\demo{Proof}
Statement (i) follows from the fact that
$\lambda$ is an eigenvalue of $A_{B}$ if and only if $\lambda^{d}$ is an
eigenvalue of $(A_{\Gamma}^{d})_{B_{d}(0)}$ and (ii) is proved in the same way.
\qed
\enddemo

\demo{Proof of Theorem 1.1}
Setting $t=0$ in (4.1), one obtains
$$logDetA^{d}-logDet(A_{\Gamma}^{d},B_{d}(0))=\tilde{c}+logDetR_{d}(0).$$
By Lemma 4.2, $log(DetA)^{d}-log(Det(A_{\Gamma},B))^{d}=\tilde{c}
+log(DetR)^{d}.$ Hence \newline
$logDetA= log(c)+logDet(A_{\Gamma},B)+logDetR$,
where $log(c)=-\frac{1}{d}\sum^{d-1}_{k=0}\pi_{0}(R(\alpha_{k}t))$. \newline
Using the result of Theorem 2.4, Theorem 1.1 follows.
\qed
\enddemo

\enddocument